\documentclass[preprint,3p]{elsarticle}

\usepackage{color}
\usepackage{amssymb}
\usepackage{amsmath}
\usepackage{appendix}
\usepackage{xcolor}
\usepackage{multicol}
\usepackage{float}

\newtheorem{Proposition}{Proposition}
\newtheorem{Definition}{Definition}

\newtheorem{Ex}{Example}

\newtheorem{proofhead}{Proof.}

\journal{}

\begin{document}

\begin{frontmatter}

\title{Modified Greenwood statistic and its application for statistical testing}

\author[inst1]{Katarzyna Skowronek}
\author[inst2]{Marek Arendarczyk}
\author[inst3]{Radosław Zimroz}
\author[inst1]{Agnieszka Wyłomańska}

\affiliation[inst1]{organization={Faculty of Pure and Applied Mathematics, Hugo Steinhaus Center, Wroclaw University of Science and Technology},
            addressline={Hoene-Wronskiego 13c}, 
            city={Wroclaw},
            postcode={50-376}, 
            country={Poland}}
                       \affiliation[inst2]{organization={Mathematical Institute, University of Wrocław},
            addressline={pl. Grunwaldzki 2/4 }, 
            city={Wrocław},
            postcode={50-384}, 
            country={Poland}}
\affiliation[inst3]{organization={Faculty of Geoengineering, Mining and Geology, Wroclaw University of Science and Technology},
            addressline={Na Grobli 15}, 
            city={Wroclaw},
            postcode={50-421}, 
            country={Poland}}
         \begin{abstract}
In this paper, we explore the modified Greenwood statistic, which, in contrast to the classical Greenwood statistic, is properly defined for random samples from any distribution. The classical Greenwood statistic, extensively examined in the existing literature, has found diverse and interesting applications across various domains. Furthermore, numerous modifications to the classical statistic have been proposed. The modified Greenwood statistic, as proposed and discussed in this paper, shares several key properties with its classical counterpart. Emphasizing its stochastic monotonicity within three broad classes of distributions - namely, generalized Pareto, $\alpha-$stable, and Student's t distributions - we advocate for the utilization of the modified Greenwood statistic in testing scenarios. Our exploration encompasses three distinct directions. In the first direction, we employ the modified Greenwood statistic for Gaussian distribution testing. Our empirical results compellingly illustrate that the proposed approach consistently outperforms alternative goodness-of-fit tests documented in the literature, particularly exhibiting superior efficacy for small sample sizes. The second considered problem involves testing the infinite-variance distribution of a given random sample.  The last proposition suggests using the modified Greenwood statistic for testing of a given distribution. The presented simulation study strongly supports the efficiency of the proposed approach in the considered problems. Theoretical results and power simulation studies are further validated by real data analysis.   

\end{abstract}

\begin{keyword}

Greenwood statistic \sep  testing Gaussianity \sep testing heavy-tailed distributions \sep Monte Carlo simulations \sep real data analysis

\end{keyword}

\end{frontmatter}

\section{Introduction}
This paper explores the modified Greenwood statistic and its use in addressing statistical testing issues. The modification made to the traditional statistic enables its application to random samples from any  distribution, whereas the original Greenwood statistic was specifically designed for random samples with positive values. The relative simple form of the introduced statistic allows for the analysis of its theoretical properties and its application in various areas of interest. 

The Greenwood statistic was introduced by M. Greenwood in 1946~\cite{medicine} and further was discussed in~\cite{Moran_1947}, where some asymptotic properties and moments of its distribution were analysed for exponentially distributed random sample. From that time many authors have analyzed this statistic and its various extensions from a theoretical perspective. For instance, Greenwood statistic and its modifications were introduced in the context of testing for the exponential and uniform distribution~(see, e.g., \cite{goodness-fit}). In~\cite{albrecher2007} and~\cite{albrecher2009} authors derived asymptotic properties of the moments of the Greenwood statistic. Also the asymptotic behaviour of the distribution of the Greenwood statistic with respect to the existence of the moments of the underlying distribution were derived in~\cite{Albrecher_2010}. In~\cite{arendarczyk2022} authors proved the stochastic monotonicity of the Greenwood statistic under the assumption of star-shaped stochastic monotonicity of an underling random sample. The Greenwood statistic was also used in~\cite{arendarczyk2023} to introduce a test for inference of the tail index in context of generalized Pareto distribution. Some of important works related to Greenwood statistic also include its applications in testing the Taylor's law, see~\cite{taylor-law1,taylor-law2}.  Recently, in~\cite{ALBRECHER2022} authors introduced a generalization of the Greenwood statistic and analyzed the asymptotic properties of the modified Greenwood statistics for regularly varying distributions. 
Due to the simple form of the Greenwood statistic and its modifications, they were applied in various real-data problems. The classical examples include applications in the analysis of clustering events either in space or time, namely in medicine and epidemiology~\cite{medicine,epidemology}, genetics and genomics~\cite{genetics,genomics}, biology~\cite{biology}, economics and insurance~\cite{economics,insurance}, hydrology~\cite{hydrology}, optimization~\cite{optimization}, physics and materials science~\cite{physics,materials}, anomaly detection~\cite{anomaly}, internet traffic monitoring~\cite{internet} and even in athletics~\cite{athletics}.

In this paper we discuss the theoretical properties of the introduced modification of the Greenwood statistic paying a particular attention to its stochastic ordering for three general classes of distributions, namely generalized Pareto distribution~\cite{emberchts1997},  $\alpha$-stable distribution~\cite{stable1}, and Student’s t distribution~\cite{student2}. Here we extended the last class, namely Student's t distribution, by the Gaussian distribution and thus, we examine the properties of the modified Greenwood statistic for such class.  The mentioned above classes of distributions cover the  light and heavy-tailed family of distributions. They are crucial in the probability theory as well as real data applications. 

In the context of application of Greenwood statistic, the generalized Pareto distribution was previously discussed in~\cite{arendarczyk2022,arendarczyk2023}. The distribution is defined through probability density function with three parameters, of which $\gamma \in \mathbb{R}$ is the most relevant parameter responsible for heaviness of the tail distribution (see, e.g.,~\cite{emberchts1997}). If $\gamma \leq 0$ the distribution has lighter tail, in particular for $\gamma=0$ the generalized Pareto distribution reduces to exponential distribution, and eventually for $\gamma \geq 0.5$ the variance of corresponding random variable does not exist. 

The $\alpha$-stable distributions are defined by four parameters, with the stability index $\alpha \in(0,2]$ considered one of the most crucial. The stability index is responsible for the heavy-tailed behavior, and the smaller values of $\alpha$ correspond to a higher probability of the associated random variable taking extreme values. When $\alpha<2$, $\alpha-$stable distributions fall within the category of heavy-tailed distributions, with the corresponding random variable exhibiting infinite variance in this scenario. On the contrary, $\alpha-$stable distributions can be viewed as an extension of the Gaussian distribution, converging to it when $\alpha=2$ (see, e.g.,~\cite{stable3}). For more details we refer the readers to classical books on $\alpha-$stable distributed signals and models, such as \cite{shao22,alek_book,non_gauss,Nolan2020}.

The class of Student's t distributions is defined through the probability density function with the parameter $\nu > 0$ (number of degrees of freedom) responsible for the tail behaviour of corresponding random variable~\cite{student2}. For $\nu > 2$, the variance of the distribution exists while it is infinite otherwise. 

As mentioned, in this paper, the modified Greenwood statistic is applied for the testing problem. Here we present three different  directions.  The first one, is related to the classical problem of testing the Gaussian distribution. In this case one can also apply the modified Greenwood statistic and propose the goodness-of-fit test. For Gaussian distribution testing, a widely used  is the Shapiro-Wilk test~\cite{SWtest} and its several extensions, see~\cite{SWmod1,SWmod2}. Another commonly utilized tests for Gaussianity are based on skewness and kurtosis, namely Jarque-Bera test~\cite{JBtest} and D'agostino-Pearson test~\cite{DPtest}. Several approaches  assessing the empirical distribution function of the random sample have been introduced to test for Gaussian distribution, e.g. Kolmogorov-Smirnov test~\cite{KStest}, Cram\'er-von-Mises test~\cite{CMtest}, Kuiper test~\cite{Kuipertest}, Watson test~\cite{Watsontest}, Anderson-Darling test~\cite{ADtest} and Lilliefors test~\cite{LFtest} (as well as its extensions~\cite{LFmod1}). In addition, several comparative studies have analyzed the effectiveness of tests for Gaussianity, e.g.~\cite{Iskander,NormalityRev2,NormalityRev1}. We also highlight the recently proposed approach based on the conditional variance statistic \cite{JelPit2018}.

The second considered problem when the modified Greenwood statistic is proposed to be applied is the testing of infinite-variance distribution for given random sample. This problem is much more general than testing a specific distribution, however it was also considered in the literature. A general test for infinite moments was introduced in~\cite{infvar1}, where the author constructed statistic that diverges if a $k$th moment is infinite and converges otherwise. Another general bootstrap-based test for finite moments was introduced in~\cite{infvar2}. Both tests can be utilized to verify whether the distribution's  variance exists. In several studies~\cite{ECEM2015,maraj2023,mssp2023} empirical cumulative even moment statistic was analyzed in context of testing for infinite variance, as the statistic diverges when random sample comes from infinite variance distribution. However, assessing the properties of the moments for various distributions with different properties is difficult, thus some alternative methods of detection of infinite moments can be applied, e.g. based on estimating the power-law behaviour of the tail of the distribution. For various distributions, the existence of the moments depend on the tail of the distribution, thus estimating the power-law index allows to infer about the existence of the variance (for estimation see~\cite{TailidxHill,Tailidx2004,Tailidx2012,Tailidx2014,Tailidx2020,Tailidx2023}). Let us note that the classical Greenwood statistic was also utilized in context of testing for tail index. Namely, in~\cite{arendarczyk2023}, the authors used the Greenwood statistic to estimate the confidence intervals for tail index parameter in case of generalized Pareto distribution.

The last application of the modified Greenwood statistic is a testing of a specific distribution with given parameters. The testing procedure is describe within three mentioned above classes. In case of non-Gaussian distributions, the goodness-of-fit tests are usually based on the empirical distributions  of the test statistic and rejection regions are obtained by Monte Carlo simulations. 
For generalized Pareto distribution this approach was extensively analyzed in~\cite{GenParetoTest1} and there are considered several goodness-of-fit tests, e.g., Kolmogorov-Smirnov, Cram\'er-von-Mises and Anderson-Darling tests, which in the classical versions were proposed for Gaussian distribution testing. Similar approach can be utilized in case of Student's t distribution. In case of the goodness-of-fit test for $\alpha$-stable distribution, the approach based on conditional moments was developed~\cite{stableTest1,maraj2023}. Moreover, in such case a likelihood ratio test have been introduced to discriminate between Gaussian distribution and $\alpha$-stable distribution with $\alpha \neq 2$ (see, e.g.,~\cite{Nolan2020}). For other goodness-of-fit test dedicated to $\alpha$-stable distributions we refer to~\cite{stableTest2}. 

The main novelty of this paper is to introduce the modified Greenwood statistic and discuss its main theoretical properties for three considered classes of distributions. Additionally, our aim is to demonstrate the usefulness of the introduced statistic for the testing problem that can be performed in three different contexts discussed above. The efficiency of the proposed testing methodology is verified for simulated random samples from three analyzed classes of distributions. Finally, the theoretical and simulation studies are supported by real data analysis from condition monitoring area. The presented results clearly confirm that the test based on modified Greenwood statistic outperforms other considered classical tests for Gaussian distribution. Specifically it is evident when dealing with small sample sizes in the class of $\alpha-$stable distributions.  Additionally, the modified Greenwood statistic serves as the powerful tool for testing the general infinite-variance distribution of given sample. 

The rest of the paper is organized as follows. In Section 2 we recall some important characteristics of the Greenwood statistic. In Section 3 we introduce the modification to Greenwood statistic and we present its main properties. In Section 4 we introduce application of the modified Greenwood statistic to testing problem. In Section 5 we present power simulation study of the proposed tests for three discussed classes of distributions. In Section 6 we show application of the proposed approach to real data case. Lastly, in Section 7 we present concluding remarks.

\section{Greenwod statistic }\label{sec2}

In this section, we present basic facts on Greenwood statistic $T_n$.
Let $X_1, X_2, \ldots, X_n$ be an independent identically
distributed (IID) random sample from a common nonnegative distribution.
Then statistic $T_n$ is defined as follows
\begin{eqnarray}    \label{Greenwood}
    T_n = \frac{\sum_{i=1}^n X_i^2}{\left(\sum_{i=1}^n X_i\right)^2}.
\end{eqnarray}
Originated by Greenwood \cite{medicine}, for testing exponentiality and
further used in a number of applications, statistic (\ref{Greenwood})
is intensively studied in recent years. Below, we present the most important
 properties of statistic $T_n$. \\

It should be noted that exact distribution of statistic $T_n$ is very difficult to obtain
and closed form expression for probability density function (PDF) or cumulative distribution function (CDF) is, in general, unknown even for underling random sample with exponential distribution.
Thus, vast literature on Greenwood statistic is devoted to the analysis of its asymptotic behavior
and approximation of percentiles.
In particular, under the assumption of finite fourth moment, $T_n$ is asymptotically normal with
mean $\frac{2}{n}$ and variance $\frac{4}{n^3}$. To be specific, let
$
    \tilde{T}_n = \sqrt{n}\left(\frac{n T_n}{2} - 1\right),
$
then
$
    \tilde{T}_n \stackrel{d}{\to} N,
$
where $N$ is a standard normal random variable with distribution $\mathcal{N}(0,1)$ (see, e.g., \cite{Moran_1947}). 
Furthermore, under the assumption of regular variation of distribution of underling random sample $X_1, X_2, \ldots, X_n$, asymptotic distribution of $T_n$ was studied in \cite{csorgo1988} and \cite{Albrecher_2010}.
In particular, in case where parameter of regular variation, belongs to the interval $(0,1)$, it was shown in \cite{Albrecher_2010} that $T_n$ converges, under proper normalisation, to $\frac{U}{V^2}$, where $U$ and $V$ are independent
random variables with $\alpha-$stable distributions (see, Theorem 2.1 and Remark 2.1 in \cite{Albrecher_2010}).
In case of regular variation parameter greater than $1$ and  proper normalisation, $T_n$ converges to a random variable with the $\alpha-$stable distribution
(see, Theorems 2.2 -- 2.5 in \cite{Albrecher_2010}).
Note that, as was pointed out in \cite{Moran_1947}, the convergence is very slow even in a case of limiting normal distribution. Thus, the usage of asymptotic distribution of $T_n$ for constructing critical regions is of limited value.
For the historical review of the most important results on approximation of Greenwood statistic see
supplementary material in \cite{arendarczyk2022}.

A vast literature is devoted to the analysis of many other interesting probabilistic
properties of statistic $T_n$.
By the definition (\ref{Greenwood}), we have $\frac{1}{n} \le T_n \le 1$, hence all moments of the Greenwood statistic are finite.
Another important characteristic of $T_n$, that is a direct consequence of (\ref{Greenwood}) is that its distribution is scale-invariant.
Note, that Greenwood statistic is closely related to some other important statistics, such as
sample coefficient of variation $CV_n = \sqrt{nT_n - 1}$, Sharpe’s ratio, Student t-statistic
$ST_n = \sqrt{\frac{n-1}{nT_n - 1}}$, and self-normalized sum $SN_n = \frac{1}{\sqrt{T_n}}$
(see, e.g., \cite{Albrecher_2010} and references therein).

Statistic $T_n$ is commonly used for testing exponentiality (see, e.g., \cite{arendarczyk2022} and the references therein). Equivalently it can be used as a test statistic for uniformity as discussed in  Section 1 of suplementary material in \cite{arendarczyk2022}.
Moreover, statistic $T_n$ and its modifications has been used in the context of testing for extreme domain of attraction
(see, e.g., \cite{Hasofer1992} and \cite{Neves2007} and the references therein).
These applications is closely related to clustering and heterogeneity detection, as it was pointed out in \cite{arendarczyk2022}.
In particular, when clustering is present $T_n$ tends to have values closer to one, while in opposite case of uniformly or even super-uniformly distributed data, that is without grouping or outliers, values of $T_n$ are closer to $\frac{1}{n}$.
We refer to Section 3 in \cite{arendarczyk2022} on the discussion about the connection between clustering, heavy tails and behavior of statistic $T_n$.
This property of Greenwood statistic, together with its stochastic monotonicity discussed below
provide justification for $T_n$ being effective tool for discriminating between light and
heavy tailed data within various classes of distributions. 

\section{Modified Greenwood statistic}

\subsection{General properties}
Due to its simple form and useful properties, as discussed in the previous sections, statistic $T_n$ is commonly used for solving numerous real-data problems. However, the main limitation of the usage of Greenwood statistic $T_n$ is that it can not be applied for general real valued samples that commonly appears in many important applications. For resolving this problem, in this section we propose a modification of Greenwood statistic called modified Greenwood statistic, defined as
\begin{eqnarray} \label{def.sn}
        S_n := \frac{\sum_{i=1}^n |X_i|^2}{\left(\sum_{i=1}^n |X_i|\right)^2}.
\end{eqnarray}
Although, proposed statistic can be used for analyzing any real valued samples, we concentrate on its value for testing within families of distributions symmetric with respect to $0$. In particular, as we show in this contribution, statistic $S_n$ exhibit high efficiency in detecting if a sample comes from $\alpha$-stable or Student's t distribution when compering with Gaussian distribution, outperforming other traditionally used test statistics.

\begin{figure} 
    \centering
    \includegraphics[width = 0.9\textwidth]{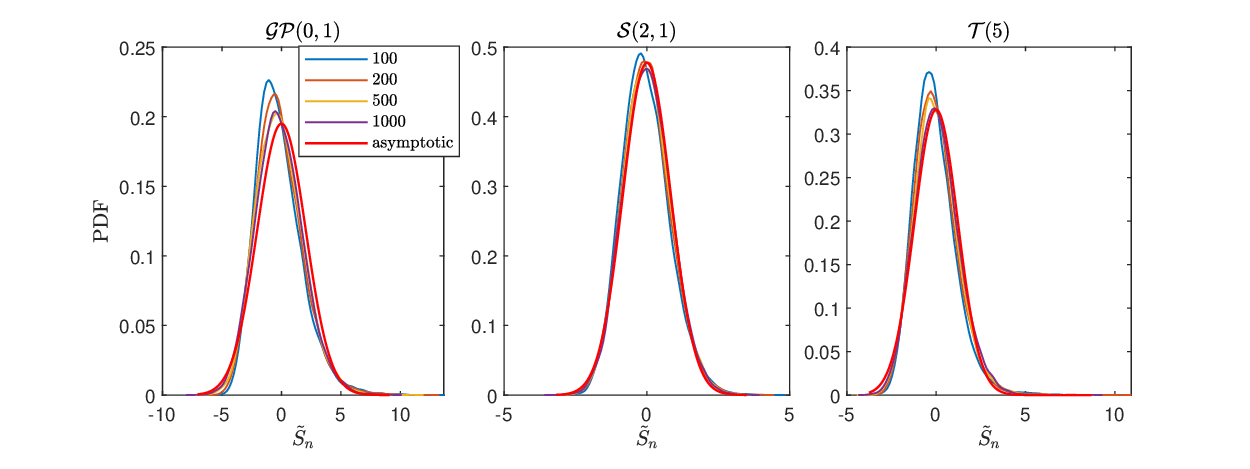}
    \caption{Comparison of PDF of asymptotic distribution of $\tilde{S}_n$ and empirical PDFs of $\tilde{S}_n$, for sample lengths $n=100,200,500,1000$ for three considered distributions. Empirical distributions of statistic $\tilde{S}_n$ are obtained based on $10000$ Monte Carlo simulations.}
    \label{fig.limit.hist}
\end{figure}

\begin{figure} 
    \centering
    \includegraphics[width = 0.9\textwidth]{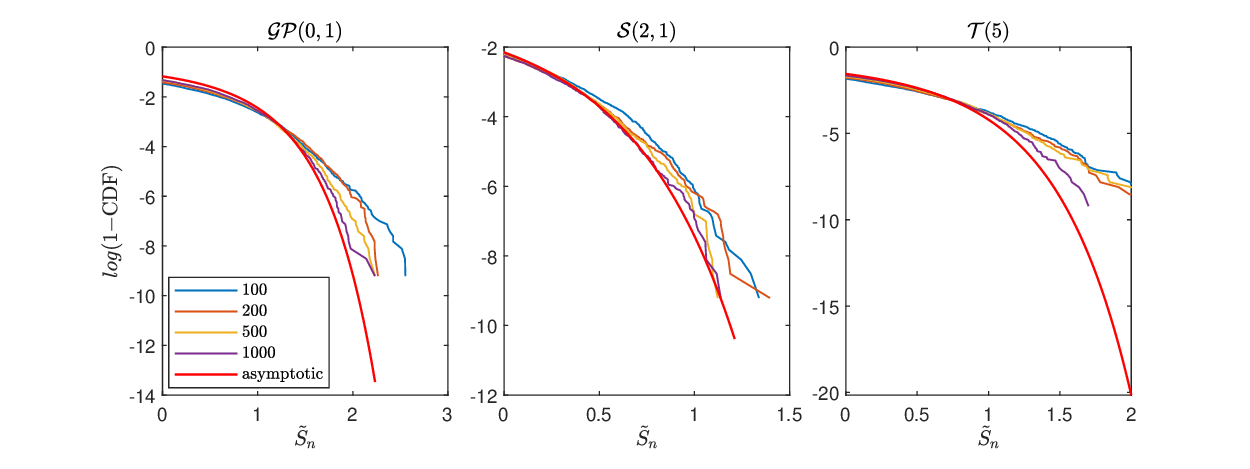}
    \caption{Comparison of tail of asymptotic distribution of $\tilde{S}_n$ and empirical tails of $\tilde{S}_n$, for sample lengths $n=100,200,500,1000$ for three considered distributions. Empirical distributions of statistic $\tilde{S}_n$ are obtained based on $10000$ Monte Carlo simulations.}
    \label{fig.limit.log}
\end{figure}

Statistic $S_n$ preserve most of the fundamental properties of Greenwood statistic. In particular, it is scale invariant.  
Moreover, note that $S_n$ is bounded, as $\frac{1}{n} \le S_n \le 1$. Consequently, all moments of $S_n$ exists even for the case where this property does not hold for underling random sample. Additionally, under assumption of finite fourth moment $S_n$ is asymptotically normal.
Note that, as was pointed out in the previous section the convergence of statistic $T_n$ is very slow. The same property persists for the modified Greenwood statistics, see Fig. \ref{fig.limit.hist} and Fig. \ref{fig.limit.log}. In Fig.\ref{fig.limit.hist}  we present the comparison of PDF of asymptotic distribution of normalized $S_n$ (denoted by $\tilde{S}_n$) and empirical PDFs  of $\tilde{S}_n$ for three considered distributions, see \ref{AppA}. Analogously, in Fig. \ref{fig.limit.log} we present the comparison of the corresponding distributions' tails (i.e., 1-CDF).  The discussed above figures clearly show that  the usage of asymptotic distribution of $S_n$ for constructing critical regions is of limited value. In the following sections we present approach for testing procedures that goes in line with the one developed in \cite{arendarczyk2022} and \cite{arendarczyk2023} and is based on stochastic monotonicity of statistics $S_n$.

\begin{figure} 
    \centering
    \includegraphics[width=0.8\textwidth]{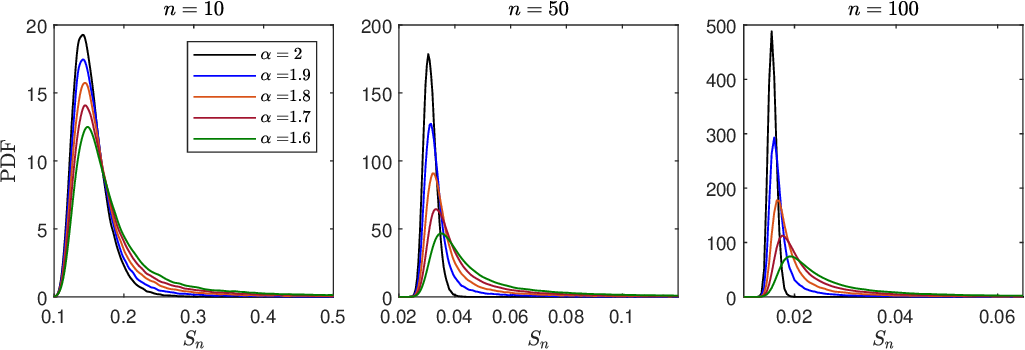}
    \caption{Comparison of the PDFs of $S_n$ for $\alpha-$stable distribution for  different values of  $\alpha$ for sample lengths $n=10, 50, 100$. The results are obtained based on $10000$ Monte Carlo simulations.}
    \label{mass.shift.fig}
\end{figure}

\subsection{Stochastic monotonicity} 
One of the important properties of statistic $S_n$ that plays a crucial role in applications
and justifies its usage as a test statistic, is its stochastic ordering.
\begin{Definition} \label{OrderDef}
Let $X$ and $Y$ be independent random variables with CDFs $F_X(\cdot)$ and $F_Y(\cdot)$, respectively. We say that $X$ is stochastically smaller than $Y$ (in standard sense), denoted by
$X \le_{\rm st} Y$, whenever, for each $x > 0$,
$
        F_Y(x) \le F_X(x).
$
\end{Definition}

It is known that stochastic ordering of Greenwood statistic $T_n$ can be obtained under assumption of star-shaped ordering of distribution of underling random sample $X_1, X_2, \ldots, X_n$. We refer the reader to Theorem 1 in \cite{arendarczyk2022} for the detailed proof of this fact.
Herein, we demonstrate that statistic $S_n$ exhibits the same stochastic ordering property. Before we present the main result of this section let us recall
two equivalent definitions of star-shaped order commonly used in the literature.

\begin{Definition} \label{StarOrderDef}
Let $X$ and $Y$ be independent random variables with CDFs $F_X(\cdot)$ and $F_Y(\cdot)$, respectively. We say that $X$ is smaller than $Y$ in the star-shaped order, denoted by
$X \le_\ast Y$, whenever
$
        g(x) := F_Y^{-1}(F_X(x))
$
is a star-shaped function, i.e., $g(\alpha x) \le \alpha g(x)$, for any, $x > 0, \alpha \in [0,1]$
(see Definition in Section 1 in \cite{Rivest_1982}).\\
Equivalently we say that $X \le_\ast Y$ if
$
        \frac{F_Y^{-1}(u)}{F_X^{-1}(u)},
$
is non-decreasing in $u \in (0,1)$, where $F_X^{-1}(\cdot)$ and $F_Y^{-1}(\cdot)$ are quantile functions of $X$ and $Y$, respectively. (see, e.g., Proposition 1 in \cite{Rivest_1982}).

\end{Definition}

Based on the assumption of star-shaped ordering of distribution of random sample $X_1, X_2, \ldots, X_n$, standard stochastic ordering of statistic $S_n$ can be shown.

\begin{Proposition} \label{SnMain}
Let $\{\mathcal{P}_\theta, \,\, \theta \in \Theta \subset \mathbb R\}$ be a family of star-shaped ordered absolutely continuous, symmetric probability distributions on $\mathbb R$, that is such that if
$\theta_1\leq \theta_2$ we have $X^{(\theta_1)}\leq_{\ast} X^{(\theta_2)}$, with  $X^{(\theta_i)}\sim \mathcal{P}_{\theta_i}$, $i=1,2$. Then, for each $n\geq 2$, we have
\begin{equation}
\label{stoch.order.tn}
S_n^{(\theta_1)}\leq_{\mbox{\rm st}}S_n^{(\theta_2)} \,\,\, \mbox{whenever $\theta_1\leq \theta_2$, $\theta_1, \theta_2\in \Theta$},
\end{equation}
where $S_n^{(\theta)}$ is given by (\ref{def.sn}) with the $\{X_i\}$ having a common distribution $\mathcal{P}_\theta$.
\end{Proposition}

The proof of Proposition \ref{SnMain} is presented in \ref{SnMain.proof}.

Moreover, observe that in case of nonnegative random variables $X_1, X_2, \ldots, X_n$ statistics $S_n$ and $T_n$ are equivalent. Thus, in this case, due to Theorem 1 in \cite{arendarczyk2022} stochastic order of $S_n$ is also preserved.

\begin{Proposition} \label{SnPositive}
Let $\{\mathcal{P}_\theta, \,\, \theta \in \Theta \subset \mathbb R\}$ be a family of star-shaped ordered absolutely continuous probability distributions on $\mathbb R_+$, that is such that if
$\theta_1\leq \theta_2$ we have $X^{(\theta_1)}\leq_{\ast} X^{(\theta_2)}$, with  $X^{(\theta_i)}\sim \mathcal{P}_{\theta_i}$, $i=1,2$. Then, for each $n\geq 2$, we have
\begin{equation}
S_n^{(\theta_1)}\leq_{\mbox{\rm st}}S_n^{(\theta_2)} \,\,\, \mbox{whenever $\theta_1\leq \theta_2$, $\theta_1, \theta_2\in \Theta$},
\end{equation}
where $S_n^{(\theta)}$ is given by (\ref{def.sn}) with the $\{X_i\}$ having a common distribution $\mathcal{P}_\theta$.
\end{Proposition}

In the subsequent part of this section we concentrate on three analyzed classes of distributions of underling random samples $X_1, X_2, \ldots, X_n$: generalized Pareto distribution $\mathcal{G}(\gamma, \delta)$, $\alpha$-stable distribution $\mathcal{S}(\alpha, \sigma)$, and Student's t distribution $\mathcal{T}(\nu)$ (see Appendix A for definitions).
Notice, that for the class $\mathcal{G}(\gamma, \delta)$ star-shaped ordering and thus standard ordering of corresponding statistic $T_n$, was proved in \cite{arendarczyk2022} and was extensively studied in the context of clustering detection and testing in \cite{arendarczyk2022} and \cite{arendarczyk2023}. We call off this result in Proposition \ref{SnOrder} (i) and refer the reader to Proposition 1 (i) in \cite{arendarczyk2022}, for its detailed proof.

\begin{Proposition} \label{SnOrder}
Let $S_n$ be a modified Greenwood statistic defined in (\ref{def.sn}), where $X_1, X_2, \ldots, X_n$ are IID random variables with a common distribution.

\begin{itemize}
\item[{\rm (i)}] If $X_i \sim \mathcal{GP}(\gamma, \delta), \gamma \in \mathbb{R}, \delta  >0$, then $S_n$ is stochastically increasing with respect to parameter $\gamma$.
\item[{\rm (ii)}] If $X_i \sim S(\alpha, \sigma), \alpha \in (0,2], \sigma > 0$, then $S_n$ is stochastically decreasing with respect to parameter $\alpha$.
\item[{\rm (iii)}] If $X_i \sim \mathcal{T}(\nu), \nu \in \mathbb{N}\cup\{\infty\}$, then $S_n$ is stochastically decreasing with respect to parameter $\nu$.
\end{itemize}
\end{Proposition}

The proof of Proposition \ref{SnOrder} is presented in \ref{SnOrder.proof}.

In Fig. \ref{mass.shift.fig} mass shift of probability is presented in the case of $\alpha$-stable distribution, as an example illustrating stochastic behavior of statistic $S_n$, demonstrated in Proposition 3 (ii). It is pointed out that the mass of probability shifts to the right as stability index moves away from 2. Similar behavior of statistic $S_n$ can be observed for generalised Pareto distribution of underlying random sample and was discussed in \cite{arendarczyk2023}.
Such behavior of the test statistic $S_n$ is strictly connected with the heaviness of the tail of distribution of underlying random sample and presence of clustering in the dataset. In particular, heavier tail in underlying random sample leads to larger values of statistic $S_n$. We refer to Section 3 in \cite{arendarczyk2023} for the broad discussion on a connection between heavy tails, clustering and Greenwood statistic. Finally, let us pointed out that stochastic monotonicity proved in Proposition \ref{SnOrder} justifies a construction of rejection regions for the tests introduced in the next section.  

\section{Application of modified Greenwood statistic for testing problem}
In this section, we demonstrate the application of the Greenwood statistic for the testing problem. We showcase the universality of the Greenwood statistic and discuss three versions of the statistical test: a test for Gaussian distribution,  a test for an infinite-variance distribution, and a test for a given distribution with a specific value of the parameter responsible for the heavy-tailed behavior. Although, the methodology presented in this section is described in a general form, in further analysis, we demonstrate its usefulness for three considered classes of distributions, namely $\alpha-$stable $\mathcal{S}(\alpha,\sigma)$, Student's t $\mathcal{T}(\nu)$  and generalized Pareto $\mathcal{GP}(\gamma,\delta)$.
\subsection{Testing of Gaussian distribution}\label{gaussian}
One of the proposed applications of Greenwood statistic is its utilization in the problem of the testing of Gaussian distribution. Let us consider the class of distributions $\mathcal{P}_\theta$, where $\theta\in \Theta$ is the distribution's parameter and $F_{\theta}(\cdot)$ denotes the CDF of the corresponding distribution. Moreover, we assume that the Gaussian distribution belongs to $\mathcal{P}_\theta$ and is characterized by $\theta^*$ parameter. Here we discuss the problem of testing Gaussian distribution versus heavy-tailed distribution. For such a case there are two possible scenarios, namely  in scenario (1) for each $\theta \in \Theta$ we have $\theta^*\leq \theta$ while in scenario (2) we have $\theta^*\geq \theta$.

In the considered problem of Gaussian distribution testing we formulate  the following $\mathcal{H}_0$ and $\mathcal{H}_1$ hypotheses, depending on the possible scenarios, respectively
\begin{eqnarray}\label{hipoteza}\mathcal{H}_0: ~\theta =\theta^{*},~~\mathcal{H}_1:~\theta^*< \theta  \\
\label{hipotezaGreater}
\mathcal{H}_0: ~\theta =\theta^{*},~~\mathcal{H}_1:~\theta^*> \theta.\end{eqnarray}
Knowing that, $S_n$ is stochastically monotone  with respect to $\theta$ within given class, we can consider it as a test statistic. For mentioned above scenarios  we consider one-sided  testing procedure meaning that the test statistic obtained for a single trajectory is compared only with one critical value and the result is the base to reject the $\mathcal{H}_0$.  The rejection region of the tests corresponding to  scenarios (1) and (2) is as follows
\begin{eqnarray} \label{MGT1.rejreg}
\left[{Q}_c(n),1\right],\end{eqnarray}
\noindent where ${Q}_p(n)$ is the theoretical quantile of order $p$ from the distribution of the test statistic under $\mathcal{H}_0$ hypothesis and $n$ is a sample length. As the distribution of $S_n$ in general is not known in the rejection region (\ref{MGT1.rejreg}) we take the empirical quantile, i.e.,  $\hat{Q}_c(n)$. To this end we simulate $M$ sample trajectories of $X_1,X_2,\ldots,X_n$ from Gaussian distribution and for each of them we calculate the value of the modified Greenwood statistic. Finally, we obtain $M$ realisations of the distribution corresponding to $S_n$. 

In further analysis, the tests corresponding to the above mentioned scenarios are denoted as $MG_1$ test and $MG_2$ test, respectively. The described above procedure can be applied for $\mathcal{S}(\alpha,\sigma)$ and $\mathcal{T}(\nu)$ distributions since in such classes the modified Greenwood statistic is stochastically monotone with respect to appropriate parameters. Additionally, $\mathcal{S}(2, \sigma)$ and $\mathcal{T}(\infty)$ distributions are Gaussian (see Appendix A).  For such two considered cases we apply $MG_2$ test. Following, in Propositions \ref{MGT1.prop} and \ref{MGT2.prop} we present fundamental properties of the test $MG_2$ applied for the distributions $\mathcal{S}(\alpha, \sigma)$ and $\mathcal{T}(\nu)$, respectively.

\begin{Proposition} \label{MGT1.prop}
The test $MG_2$ for the hypotheses (\ref{hipotezaGreater}) applied for the distribution 
$\mathcal{S}(\alpha, \sigma)$ has size $c$, is unbiased and has decreasing power functions with respect to parameter $\alpha$.
\end{Proposition}

The proof of Proposition \ref{MGT1.prop} is presented in \ref{MGT1.prop.proof}.

\begin{Proposition} \label{MGT2.prop}
The test $MG_2$ for the hypotheses (\ref{hipotezaGreater}) applied for the distribution 
$\mathcal{T}(\nu)$ has size $c$, is unbiased and has increasing power functions with respect to parameter $\nu$.
\end{Proposition}

The proof of Proposition \ref{MGT2.prop} shall be omitted as it is similar to the proof of Proposition \ref{MGT1.prop}.

\subsection{Testing of infinite-variance distribution}\label{vartest}

In this part, we present the application of the modified Greenwood statistic in the procedure for testing of infinite-variance distribution which is much more general problem than testing of Gaussian distribution. As in previous case, let us assume that random sample comes from the family of distributions  with CDF $F_\theta( \cdot )$. We assume the statistic $S_n$ is stochastically monotone within this class with respect to $\theta \in \Theta$. Here we assume that $\theta$ is the distribution's parameter responsible for finiteness of variance. The hypotheses of the test depend on the behavior of the tail of the distribution with respect to $\theta$. Let us assume that $\theta^* \in \Theta$ is the boundary between finite-variance and infinite-variance distribution. Here we consider two tests corresponding to two possible scenarios. In the scenario (3) the infinite-variance distribution corresponds to $\theta\geq \theta^*$ while in scenario (4) to $\theta\leq \theta^*$. Thus, for such two tests we have the following $\mathcal{H}_0$ and $\mathcal{H}_1$ hypotheses, respectively
   \begin{eqnarray} \label{hyp.var.greater}
   \mathcal{H}_0:~~\theta  \geq \theta^*,~~~\mathcal{H}_1:~\theta  < \theta^*, \\
   \label{hyp.var.less}
   \mathcal{H}_0:~\theta  \leq \theta^* ,~~~\mathcal{H}_1:~\theta  > \theta^*. \end{eqnarray}
 \noindent In this class of problems, we utilize one-sided testing procedure. Thus, the statistic obtained for a single trajectory is compared only with one critical value.  For discussed tests, the rejections regions are defined as follows
 \begin{eqnarray}\label{qq1}
     \left[\frac{1}{n}, {Q}_c(n)\right],
 \end{eqnarray}

\noindent where, similar as in the previous case,  ${Q}_p(n)$ is the theoretical $p$-quantile from the distribution of $S_n$ under $\mathcal{H}_0$ and $n$ is a sample length. Here we also apply the empirical quantile in (\ref{qq1}), i.e., $\hat{Q}_c(n)$, as the theoretical distribution of $S_n$ is unknown. In the further analysis the tests are denoted as $MG_3$ and $MG_4$, respectively. Similar as in previous case, the rejection regions with certain confidence level $c$ are constructed based on $M$ simulated trajectories of $X_1,X_2, ..., X_n$ from the distribution with CDF $F_{\theta^*} (\cdot)$.\\
\indent It is worth noting that similar approach was introduced in~\cite{arendarczyk2023} in the context of generalized Pareto distribution. In the current paper, we extend the testing procedure to Student's t distribution. Test $MG_3$ corresponds to the test (4.3) from~\cite{arendarczyk2023} in the generalized Pareto distribution class. In this paper, testing the infinite-variance distribution for Student's t class corresponds to test $MG_4$. The respective values of $\theta^*$ parameters used to calculate rejection regions are $\gamma=0.5$  and $\nu=2$  for $\mathcal{GP}(\gamma,\delta)$ and $\mathcal{T}(\nu)$ distributions, respectively.

In Proposition \ref{MGT4.prop} we present fundamental properties of the test $MG_4$ applied for the distributions
$\mathcal{T}(\nu)$.

\begin{Proposition} \label{MGT4.prop}
The test $MGT_4$ for the hypotheses (\ref{hyp.var.less}) applied for the distribution
$\mathcal{T}(\nu)$ has size $c$, is unbiased and has increasing power functions with respect to parameter $\nu$.
\end{Proposition}

The proof of Proposition \ref{MGT4.prop} shall be omitted as it is similar to the proof of Proposition \ref{MGT1.prop}.

{
Finally, let us note, that we can not test infinite-variance distribution of $\alpha-$stable class, strictly in the framework presented in the current section, however, in practical applications one may consider the inverse problem, i.e., testing of finite-variance distribution, which for $\alpha-$stable class of distributions reduces to the testing of Gaussian distribution (see Section \ref{gaussian}).} 

\subsection{Testing of given distribution with specific value of the parameter responsible for heavy-tailed behavior}\label{general}
In this case, similar as previous cases, we assume 
that considered random sample comes from the distribution with the CDF $F_{\theta}(\cdot)$, where $\theta\in \Theta$ and within this class the modified Greenwood statistic $S_n$ is stochastically monotone with respect to parameter $\theta$. In this case the $\mathcal{H}_0$ and $\mathcal{H}_1$ are defined as 
\begin{eqnarray}~~\mathcal{H}_0: ~\theta =\theta^{*},~~~\mathcal{H}_1:~\theta^*\neq \theta.
\end{eqnarray}
Here $\theta^{*}$ may correspond to any case of the considered family of distributions in contrast to the case presented in Section \ref{gaussian}, where $\theta^{*}$ was the parameter corresponding to Gaussian distribution. Additionally, we assume there exist $\theta_1,\theta_2\in \Theta$ such that $\theta_1<\theta^*<\theta_2$. 

Here we apply two-sided testing procedure. For single sample trajectory we reject $\mathcal{H}_0$ hypothesis if the test statistic is extreme, either larger than an upper critical value or smaller than a lower critical value with a given significance level $c$. To  construct the rejection region of the considered test at given confidence level $c$ based on the statistic $S_n$, we proceed similarly as in the previous cases.  We simulate $M$ sample trajectories of $X_1,X_2,\ldots,X_n$ from distribution with CDF $F_{\theta^{*}}(\cdot)$. For each simulated sample we calculate the value of the modified Greenwood statistic. Finally, we calculate the rejection region
\begin{eqnarray}
\left[\frac{1}{n}, {Q}_{c/2}(n)\right]\cup\left[{Q}_{1-c/2}(n),1\right],
\end{eqnarray}
\noindent where, ${Q}_p(n)$ is the theoretical $p$-quantile from the distribution of $S_n$ under $\mathcal{H}_0$ and $n$ is a sample length.

The described above test can be applied for all three distributions considered in this paper, namely $\mathcal{GP}(\gamma,\delta)$, $\mathcal{S}(\alpha,\sigma)$ and $\mathcal{T}(\nu)$ distributions. In those cases the parameters $\theta$ used in the general description are $\alpha, \gamma $ and $\nu$, respectively.

\section{Power simulation study}

In this section, we present results reflecting the versatility of utilizing the proposed statistic in various statistical tests. At first, we apply the test for Gaussianity described in Section ~\ref{gaussian}. Eventually, we show the efficiency of the test for infinite variance introduced in Section~\ref{vartest}. Moreover, for each test we compare the obtained results with the results of testing methods known in the literature.

\subsection{Testing of Gaussian distribution}

\begin{figure}
    \centering
    \includegraphics[width = 0.9\textwidth]{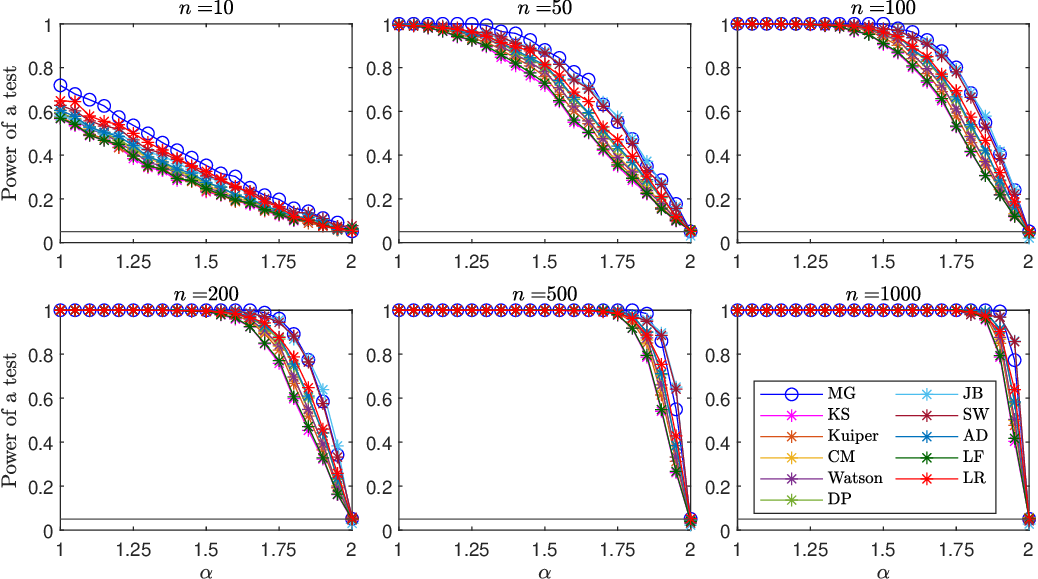}
    \caption{The power of a MG test for Gaussian distribution obtained for the class of $\alpha$-stable distributions. In this case, $\mathcal{H}_0$ corresponds to Gaussian distribution and $\mathcal{H}_1$ to $\alpha$-stable distribution with $\alpha \neq 2$. The power of MG test is compared with the power of KS, Kuiper, CM, Watson, DP, JB, SW, AD, LF and LR tests. The simulations were conducted for $\alpha \in \{1,1.05,\ldots, 2\}$ and for sample lengths $n=10,50,100,200,500,1000$. Black solid line represent the significance level $c=0.05$ selected for all analyzed tests. The results are obtained based on $2000$ Monte Carlo simulations.}
    \label{fig:stable_gauss_test}
\end{figure}

\begin{figure}[htp!]
    \centering
    \includegraphics[width = 0.9\textwidth]{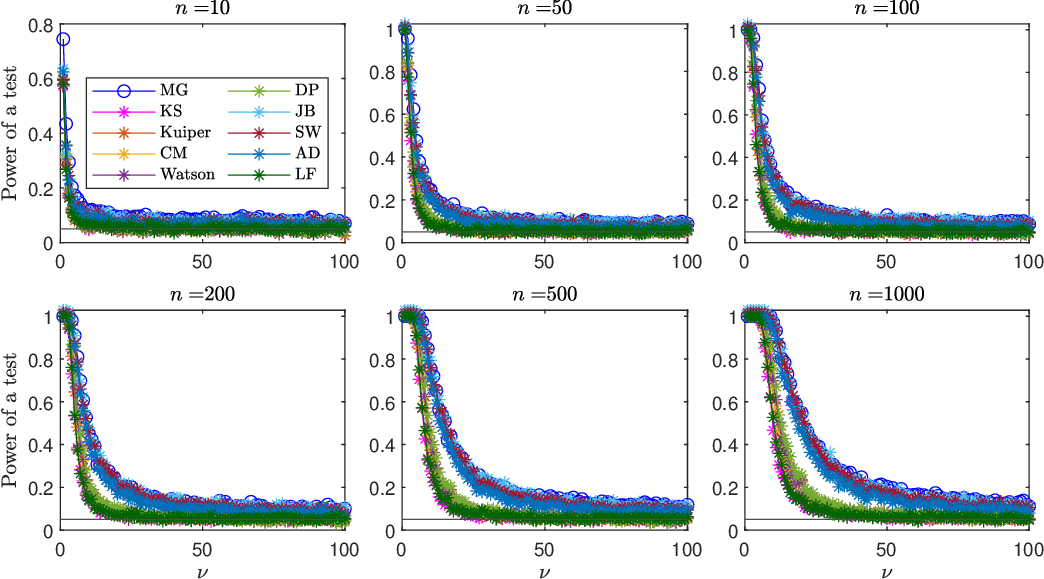}
    \caption{The power of a MG test for Gaussian distribution obtained for the class of Student's t distributions. In this case, $\mathcal{H}_0$ corresponds to Gaussian distribution and $\mathcal{H}_1$ to Student's t distribution. The power of MG test is compared with the power of KS, Kuiper, CM, Watson, DP, JB, SW, AD and LF tests. The simulations were conducted for $\nu \in \{1,2,\ldots, 100\}$ and for sample lengths $n=10,50,100,200,500,1000$. Black solid line represent the significance level $c=0.05$ selected for all analyzed tests. The results are obtained based on $2000$ Monte Carlo simulations.}
    \label{fig:tstud_gauss_test}
\end{figure}

In this part, we present results obtained for testing of the Gaussian distribution. This part is related to the test introduced in Section~\ref{gaussian}. In this case, $\mathcal{H}_0$ corresponds to Gaussian distribution of a random sample and $\mathcal{H}_1$ corresponds to non-Gaussian distribution. The testing procedure and quantiles are based on the distribution of $S_n$ obtained in $100000$ Monte Carlo simulations, thus as the number of simulations is large, the empirical distribution of $S_n$ is close to the theoretical one. The quantiles obtained in testing procedure are presented in Table~\ref{Tab:gauss} (see~\ref{tables}). Later, we denote modified Greenwood statistic-based test introduced in this paper as MG test. Moreover, we compare the results of the proposed test for Gaussianity with the tests known in the literature. We selected tests thoroughly investigated in~\cite{Iskander} and commonly used in various applications, namely Kolmogorov-Smirnov (KS) test~\cite{KStest}, Kuiper test~\cite{Kuipertest}, Cram\'er-von-Mises (CM) test~\cite{CMtest}, Watson test~\cite{Watsontest}, D'agostino Pearson (DP) test~\cite{DPtest}, Anderson-Darling (AD) test~\cite{ADtest}, Shapiro-Wilk (SW) test~\cite{SWtest} and Jarque-Bera (JB) test~\cite{JBtest}. In addition, we also compare the results with the Lillefors (LF) test~\cite{LFtest}, and in the class of the $\alpha$-stable distribution with the likelihood ratio (LR) test~\cite{Nolan2020}. 
For all of the tests, we set significance level to $c=0.05$. To prove the efficiency of the proposed test, we show the power of the tests based on 2000 Monte Carlo simulations. For $\alpha$-stable distributions, we tested the procedure for samples with stability index $\alpha \in \{1,1.05,\ldots, 2\}$. In case of Student's t distribution, we selected number of degrees of freedom in range $\nu \in \{1,2, \ldots, 100\}$. We simulated trajectories with sample lengths $n =10,50,100,200,500,1000$. The results are presented in Fig.~\ref{fig:stable_gauss_test} and~ Fig. \ref{fig:tstud_gauss_test}. \\
\indent For $\alpha$-stable distribution (see Fig.~\ref{fig:stable_gauss_test}), the power of a test introduced in this paper decreases with increasing value of stability index $\alpha$, that is as the $\alpha$-stable distribution tends to Gaussian distribution. Thus, the power increases as stability index $\alpha$ moves away from $2$, that represents Gaussian distribution. Moreover, for smaller sample sizes, the proposed MG test also performs better or at least as good as any alternative test considered in this paper. The MG test is more efficient especially for $n = 10,50$, which means that when $\alpha<2$, the MG test tends to rightfully reject $\mathcal{H}_0$ more often than other analyzed techniques. For $n=100$ and $\alpha \geq 1.85$, the JB test outperforms MG test. However, for $n=100$ and $\alpha < 1.85$, the MG test performs better than other considered tests. For $n = 200,500,1000$ the JB and SW tests slightly outperform MG test for large values of $\alpha$, however in all the considered cases, the MG test is more efficient than any other test known from the literature. 

The results obtained for the class of Student's t distributions are presented in Fig.~\ref{fig:tstud_gauss_test}. The power of MG test decreases as $\nu$ decreases, which is expected result, since the Student's t distribution tends to Gaussian distribution as $\nu\rightarrow\infty$. The proposed test performs as good as the JB, SW or AD tests, and outperforms other tests considered in this paper. Thus, as the power of MG test is large for large values of $\nu$, the test can be utilized to distinguish between near-Gaussian Student's t distributions (namely when $\nu$ is large) and Gaussian distributions.

\subsection{Testing of infinite-variance distribution}

\begin{figure}
    \centering
    \includegraphics[width = 0.9\textwidth]{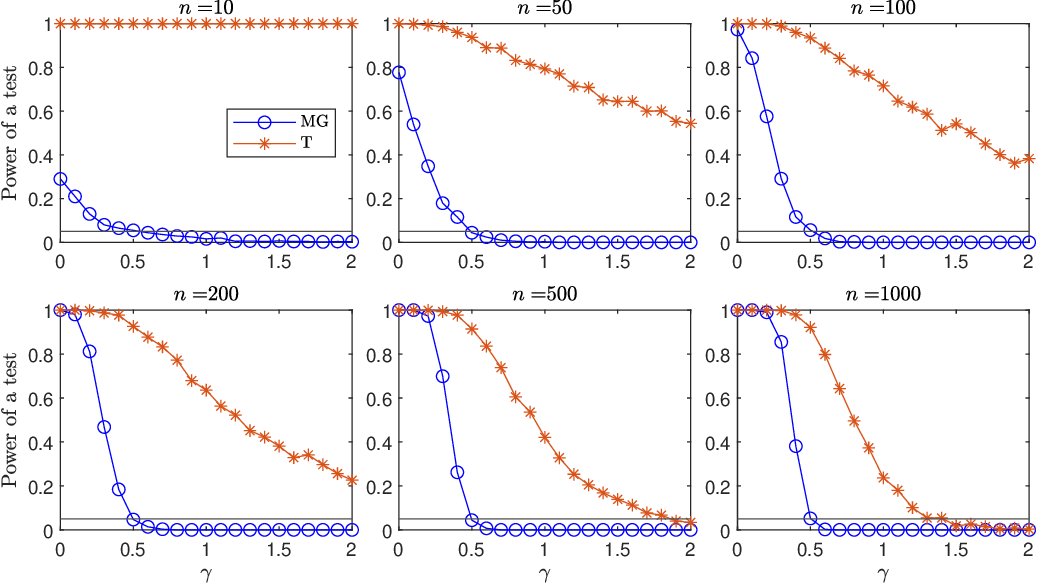}
    \caption{Results of MG test for infinite-variance for generalized Pareto distribution. In this case, $\mathcal{H}_0$ corresponds to infinite-variance distribution and $\mathcal{H}_1$ corresponds to finite-variance distribution. The power of MG test is compared with the power of T test. Power of a test was calculated based on 2000 Monte Carlo simulations. The simulations were conducted for $\gamma \in \{0,0.1,\ldots,2\}$ for sample lengths $n=10,50,100,200,500,1000$. Black solid line represents significance level $c=0.05$.}
    \label{fig:pareto-infvar}
\end{figure}

In this section, we present results obtained for testing for infinite-variance distribution for the class of generalized Pareto and Student's t distributions. In this case, the $\mathcal{H}_0$ corresponds to infinite-variance distribution, and $\mathcal{H}_1$ corresponds to finite-variance distribution. For each distribution, the quantiles of the distribution of the modified  Greenwood statistic used in testing procedures were obtained based on $100000$ Monte Carlo simulations. Quantile of $S_n$ statistic for generalized Pareto distribution was calculated for $\gamma= 0.5$, and for Student's t for $\nu=2$. The quantiles obtained in testing procedure are presented in Table~\ref{Tab:infvarGP} and~\ref{Tab:infvarT} (see Appendix~\ref{tables}). Power of a test was calculated based on $2000$ Monte Carlo simulations. Moreover, we compare the performance of Greenwood statistic based test with test known in the literature, later denoted as T test. In~\cite{infvar1}, the author proposed a statistical test to identify if the random sample comes from the distribution with infinite variance, and the test statistic of the T test under $\mathcal{H}_0$ has known theoretical distribution. For MG test and T test the significance level was set to $c=0.05$. The results of the power of the tests were obtained for the parameters responsible for the tail of the distributions: $\gamma \in \{0,0.1,\ldots,2\}$ for generalized Pareto distribution, and $\nu \in \{1,2,\ldots, 15\}$ for Student's t distribution. Simulated trajectories had lengths $n=10,50,100,200,500,1000$. \\ 
\indent For generalized Pareto distribution when $\gamma \geq 0.5$ the variance is infinite, while it is finite otherwise. In Fig.~\ref{fig:pareto-infvar}, we present results obtained for MG test and T test for infinite variance in case of generalized Pareto distribution. The power of introduced MG test decreases with increasing $\gamma$, as expected. The test outperforms the T test, even in case of smaller sample sizes. Namely, the MG test correctly does not reject $\mathcal{H}_0$ for $\gamma \geq 0.5$, and for $\gamma < 0.5$ the test rejects $\mathcal{H}_0$ in favor of $\mathcal{H}_1$. Moreover, when $n=10$ the T test rejects $\mathcal{H}_0$ for all values of $\gamma$ and fails to distinguish between infinite- and finite-variance distributed samples.

\begin{figure}[H]
    \centering
    \includegraphics[width = 0.9\textwidth]{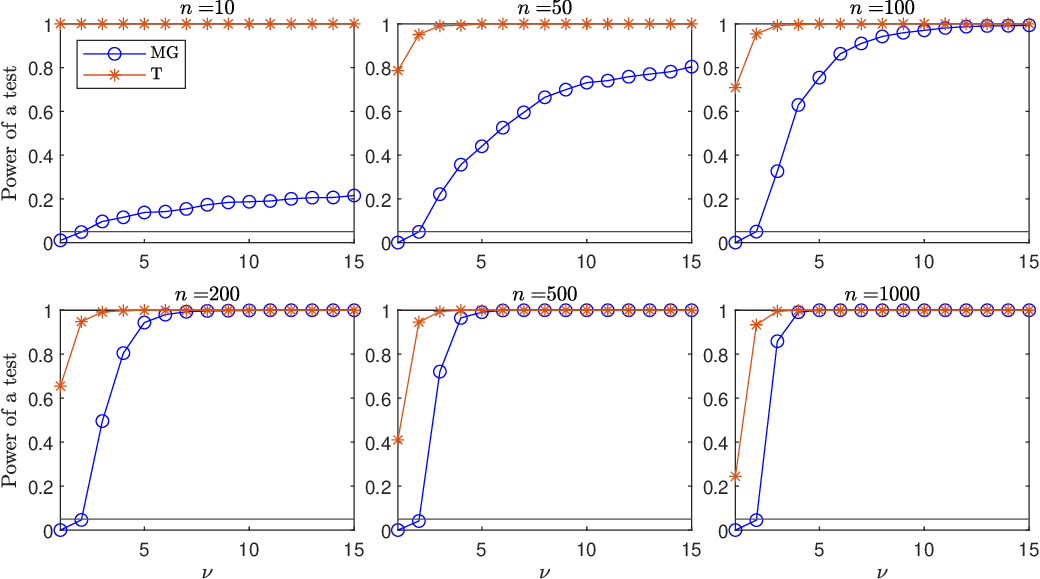}
    \caption{Results of MG test for infinite-variance for Student's t distribution. In this case, $\mathcal{H}_0$ corresponds to infinite-variance distribution and $\mathcal{H}_1$ corresponds to finite-variance distribution. The power of MG test is compared with the power of T test. Power of a test was calculated based on 2000 Monte Carlo simulations. The simulations were conducted for $\nu \in \{1,2,\ldots,15\}$ for sample lengths $n=10,50,100,200,500,1000$. Black solid line represents significance level $c=0.05$.}
    \label{fig:tstudent-infvar}
\end{figure}

\indent The power of MG and T tests for Student's t distribution are presented in Fig~\ref{fig:tstudent-infvar}. In this case for $\nu>2$ the variance exists while for $\nu \leq 2$ it is infinite. The power of proposed test increases as $\nu$ increases, as expected. As the sample size increases, the proposed in this paper test becomes more accurate, which means that for finite-variance distribution with $\nu>2$, it tends to more often rightfully reject $\mathcal{H}_0$ hypothesis in favor of $\mathcal{H}_1$. Comparing with T test, the MG test is more effective for all sample sizes, as it does not reject $\mathcal{H}_0$ for $\nu \in \{1,2\}$ when the distribution has infinite variance. Moreover, T test wrongly rejects $\mathcal{H}_0$ when $n=10,50$ for all $\nu$, failing to detect infinite-variance for small sample sizes as opposed to MG test, which does not reject $\mathcal{H}_0$ for $\nu \in \{1,2\}$ and rejects $\mathcal{H}_0$ in favor of $\mathcal{H}_1$ otherwise.

\section{Real data analysis}
 
In this section we investigate the efficiency of the proposed MG test for real data. The considered problem is related to condition monitoring in mining industry and identification of the properties of the background noise in the vibration signal. Here we apply the modified Greenwood statistic-based test for Gaussian distribution and for infinite-variance distribution of the examined random sample. 
This approach is especially useful in the analysis connected to condition monitoring, when one have to select proper tool for fault detection based on the vibration signal, see \cite{mssp2023} for more details about the examined problem from the condition monitoring perspective. The selected dataset was analyzed in~\cite{ecem2023,mssp2023}. The signal was obtained from hammer crusher in good condition. Hammer crusher is a kind of machine that is used for fragmentation of lumps of cooper ore.  The dataset consists of 25500000 observations collected during $1.7$ minutes. \\
\indent In our analysis, we apply the testing procedure to the raw signal (time domain) and to its spectrogram representation (time-frequency domain), as for the vibration signal the time-frequency representation of the data is usually examined. As the signal in the time domain has 25500000 observations, we segmented the data into $25500$ non-overlapping sub-signals, each consisting of $1000$ observations. In time-domain for such sub-signals the proposed approach is applied. 

Let us note that for lower frequency range, usually there are some periodic components~\cite{mssp2023}. Thus, in order to analyze the signal in time-frequency domain, first we applied high-pass filtering to data in time domain and analyzed higher range of frequencies from $1kHz$ to $12.5kHz$. The spectrogram $S(f,t)$ is defined as square of the short time Fourier transform of given signal~\cite{mssp2023}, namely

\begin{equation}
    S(f,t) = \left| \sum_{m=1}^n x_m w(t - m)e^{-i2\pi f \frac{m}{n}} \right|^2,
\label{eq:S}
\end{equation}
\noindent where $x_1,x_2, \dots, x_n$ are the observations in the time domain, $w(\cdot)$ is a window function, $t \in T$ is a time point and $f \in F$ is the frequency. Hence, the spectrogram is two dimensional map. It is worth noting that the number of time points $\Tilde{n}$ in the time-frequency domain depends on the windowing function, overlapping parameter and the number of the observations in the time domain. In the testing procedure vectors $S(f,t_1), S(f,t_2), \ldots, S(f,t_{\Tilde{n}})$ are considered as the random samples. Thus, for each frequency we obtain a sub-signal of length $\Tilde{n}$. As the MG test is valid for independent data, in our analysis the overlapping parameter is set to 0. Selected windowing function was $kaiser(2000,5)$ in MATLAB software. Thus, in time-frequency domain we obtained 943 sub-signals, each consisting of 1275 observations. \\
\indent For such data, we first we apply the test for Gaussian distribution and then for the infinite-variance distribution in order to confirm the results presented in \cite{ecem2023,mssp2023}. In case of modified Greenwood statistic-based test for Gaussian distribution, $\mathcal{H}_0$ corresponds to Gaussian distribution and $\mathcal{H}_1$ corresponds to non-Gaussian distribution. The distribution of $S_n$ statistic under $\mathcal{H}_0$ hypothesis is obtained based on $10000$ Monte Carlo simulations. In the MG test for infinite variance, $\mathcal{H}_0$ corresponds to infinite-variance distribution and $\mathcal{H}_1$ corresponds to finite-variance distribution. We apply the testing procedure for two classes of distributions: generalized Pareto and Student's t distributions. The distribution of $S_n$ statistic under $\mathcal{H}_0$ hypothesis is obtained based on $10000$ Monte Carlo simulations for $\gamma=0.5$ and $\nu=2$ for generalized Pareto and Student's t distribution, respectively. In case of analysis in the time-frequency domain, the distribution of $S_n$ was obtained for spectrograms of simulated random samples from Gaussian distribution, and for $\gamma=0.5$ and $\nu=2$ for generalized Pareto and Student's t distribution, respectively. \\ 
\indent The results of the testing procedures are presented in Table~\ref{Tab:realdata2} (for Gaussian distribution testing) and Table~\ref{Tab:realdata} (for infinite-variance distribution testing). In case of testing for Gaussian distribution, the MG test rejects $\mathcal{H}_0$ for $89\%$ sub-signals in the time domain and for all of the frequencies in the time-frequency domain. Thus, we conclude that the signal is not Gaussian-distributed. In case of MG test for infinite variance, under the assumption of generalized Pareto distribution, in time domain the $\mathcal{H}_0$ was rejected for $45\%$ of the sub-signals and in time-frequency domain $\mathcal{H}_0$ was rejected for $38\%$ frequencies. Thus, in time domain for the MG test with rejection region of $S_n$ calculated for generalized Pareto distribution, the results are ambiguous. For MG test based on distribution of $S_n$ statistic obtained for Student's t distribution, in time domain the $\mathcal{H}_0$ was rejected for $25\%$ of the sub-signals, and in spectrogram representation $\mathcal{H}_0$ was rejected for $11\%$ of frequencies. Thus, the distribution of the noise of the signal can be classified as the infinite-variance as for both time and time-frequency domain the majority of the sub-signals is classified as infinite-variance distributed. The obtained results confirm the results from both~\cite{ecem2023} and~\cite{mssp2023}, where the signal also was classified as infinite-variance distributed.

\begin{table}[htp!]
\centering
\caption{Results obtained in testing for Gaussian distribution in time domain and time-frequency domain. In this case, $\mathcal{H}_0$ corresponds to Gaussian distribution, $\mathcal{H}_1$ corresponds to non-Gaussian distribution. We present the percentage of sub-signals for which the $\mathcal{H}_0$ was rejected. }
\begin{tabular}{lll}
                    & MG &  \\ \cline{1-3} 
[\%] time domain            & 89                &   \\ \cline{1-3}
[\%] time-frequency domain               & 100                &   \\
                
\end{tabular}
\label{Tab:realdata2}
\end{table}

\begin{table}[htp!]
\centering
\caption{Results obtained in testing for infinite variance in time domain and time-frequency domain. In this case, $\mathcal{H}_0$ corresponds to infinite-variance distribution, $\mathcal{H}_1$ corresponds to finite-variance distribution. We present the percentage of sub-signals for which the $\mathcal{H}_0$ was rejected.
}
\begin{tabular}{llll}
                    & MG$_{\mathcal{GP}}$ & MG$_{\mathcal{T}}$ & \\ \cline{1-4} 
[\%] time domain            & 45                &  25               &  \\ \cline{1-4}
[\%] time-frequency domain               & 38                & 11               &  \\
                
\end{tabular}
\label{Tab:realdata}
\end{table}

\begin{figure}[H]
    \centering
    \includegraphics[width=0.7\textwidth]{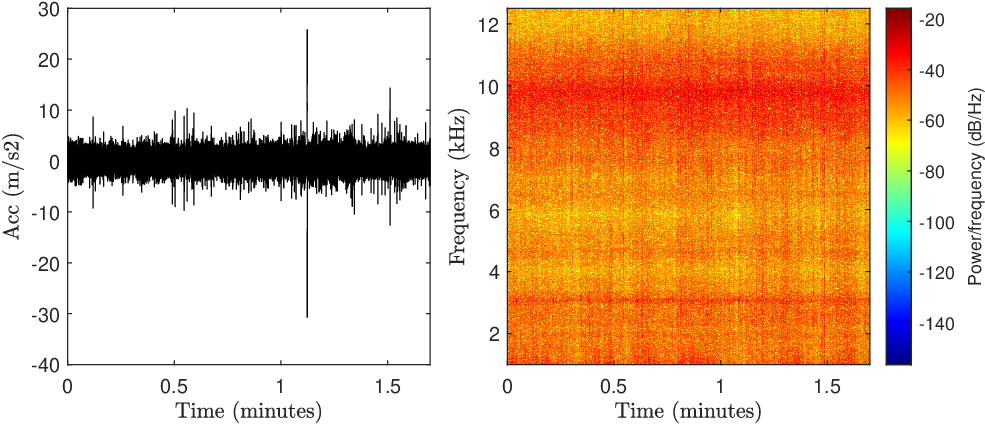}
    \caption{Analyzed signal in time and time-frequency domain. Left panel presents signal obtained from the crusher (time domain). Right panel presents the spectrogram of the analyzed signal (time-frequency domain).}
    \label{fig:condition}
\end{figure}

\section{Summary and Conclusions}
In this paper, we discuss the modified Greenwood statistic, which is a simple extension of the popular Greenwood statistic widely discussed in the literature. The modified Greenwood statistic can be applied to random samples from any distribution, thus its applicability is much wider than in the case of the classical Greenwood statistic, which is defined only for samples with positive distributions. Here, we discuss the main properties of the proposed statistic for general random samples, with particular attention paid to the stochastic monotonicity within three general classes of distributions: namely, the generalized Pareto, $\alpha$-stable, and Student's t distributions. The class of Student's t distributions is enriched with the Gaussian distribution, which is considered in this class as the limiting distribution. The stochastic monotonicity of the modified Greenwood statistics is discussed with respect to the parameters responsible for the tail behavior in the considered classes of distributions, namely $\gamma$ (in the case of the generalized Pareto distribution), the stability index $\alpha$ (in the case of the $\alpha$-stable distribution), and the number of degrees of freedom $\nu$ (in the case of the Student's t distribution).

The proven theoretical properties of the modified Greenwood statistic within three general classes of distributions enable us to propose it as a test statistic. In this paper, we describe three testing scenarios where the introduced statistic can be applied. In the first scenario, the modified Greenwood statistic is applied in the goodness-of-fit test for the Gaussian distribution. This is a classical problem discussed widely in the literature, with many proposed solutions. Our proposition extends the existing literature. The presented simulation study, where we analyze the power of the test for the Gaussian distribution within $\alpha$-stable, and Student's t classes of distributions, clearly confirms the superiority of the modified Greenwood statistic-based test over classical tests, especially for small sample sizes.

In the second scenario, we propose using the modified Greenwood statistic to test for the infinite-variance distribution of a given sample. This problem is much more general than testing for a specific distribution. We note that knowledge about the finiteness of the theoretical variance is crucial for selecting appropriate tools for data analysis. This point has been extensively discussed in the context of condition monitoring, as seen in \cite{mssp2023}, however, it is also crucial in many other areas of interest. Here, the modified Greenwood statistic proves to be a perfect tool for the considered problem. Through presented simulation studies for two classes of distributions, namely, generalized Pareto and Student's t, we demonstrate the efficacy of the proposed approach compared to the known method proposed in \cite{infvar1}.

The last scenario proposed in this paper involves using the modified Greenwood statistic for testing a given distribution with assumed parameters. In the simulation study, this problem is discussed in the context of testing the Gaussian distribution, so we did not perform specific simulations for the third scenario. However, we highlight the usefulness of the modified Greenwood statistic in this area as well.

To demonstrate the usefulness of the proposed approach, we consider dataset from condition monitoring area. It has been previously examined in the literature, where its specific properties were discussed. We analyzed a vibration signal to identify possible infinite-variance distributions within the data. Vibration signals are typically analyzed in time-frequency domains, often using spectrograms, and our analysis was conducted accordingly. The signal under examination was also discussed in \cite{mssp2023}, where the authors proposed a simple technique to identify non-Gaussian heavy-tailed behavior. In this paper, employing the modified Greenwood statistic, we corroborated previous findings and determined that such data correspond to an infinite-variance distribution in both the time and time-frequency domains. This result is pivotal for further analysis of such data, particularly in the context of local damage detection; for more information, refer to \cite{mssp2023}.

The simulation studies and real data analysis presented clearly confirm the usefulness of the modified Greenwood statistic in testing problems and its superiority over classical approaches.

 \section*{Acknowledgements}
The work of KS, RZ, and AW is supported by National Center of Science under Sheng2 project No. UMO-2021/40/Q/ST8/00024 "NonGauMech - New methods of processing non-stationary signals (identification, segmentation, extraction, modeling) with non-Gaussian characteristics for the purpose of monitoring complex mechanical structures". 

\section*{Conflicts of interest}
The authors declare no conflict of interest.

 \bibliographystyle{elsarticle-num} 
 \bibliography{mybibliography}

\begin{thebibliography}{10}
\expandafter\ifx\csname url\endcsname\relax
  \def\url#1{\texttt{#1}}\fi
\expandafter\ifx\csname urlprefix\endcsname\relax\def\urlprefix{URL }\fi
\expandafter\ifx\csname href\endcsname\relax
  \def\href#1#2{#2} \def\path#1{#1}\fi

\bibitem{medicine}
M.~Greenwood, The statistical study of infectious diseases, Journal of the Royal Statistical Society 109~(2) (1946) 85--110.

\bibitem{Moran_1947}
P.~A.~P. Moran, The random division of an interval, Supplement to the Journal of the Royal Statistical Society 9~(1) (1947) 92.

\bibitem{goodness-fit}
R.~B. D'Agostino, M.~A. Stephens, Goodness-of-fit techniques, Marcel Dekker, Inc., USA, 1986.

\bibitem{albrecher2007}
H.~Albrecher, J.~Teugels, Asymptotic analysis of a measure of variation, Theory of Probability and Mathematical Statistics 74 (2007) 1--10.

\bibitem{albrecher2009}
H.~Albrecher, J.~L. Teugels, K.~Scheicher, A combinatorial identity for a problem in asymptotic statistics, Applicable Analysis and Discrete Mathematics 3~(1) (2009) 64--68.

\bibitem{Albrecher_2010}
H.~Albrecher, S.~A. Ladoucette, J.~L. Teugels, Asymptotics of the sample coefficient of variation and the sample dispersion, Journal of Statistical Planning and Inference 140~(2) (2010) 358–368.

\bibitem{arendarczyk2022}
M.~Arendarczyk, T.~J. Kozubowski, A.~K. Panorska, {The Greenwood statistic, stochastic dominance, clustering and heavy tails}, Scandinavian Journal of Statistics 49~(1) (2022) 331--352.

\bibitem{arendarczyk2023}
M.~Arendarczyk, T.~J. Kozubowski, A.~K. Panorska, {A Computational Approach to Confidence Intervals and Testing for Generalized Pareto Index Using the Greenwood Statistic}, REVSTAT-Statistical Journal 21~(3) (2023) 367--388.

\bibitem{taylor-law1}
M.~BROWN, J.~E. COHEN, V.~H. D.~L. PEÑA, Taylor's law, via ratios, for some distributions with infinite mean, Journal of Applied Probability 54~(3) (2017) 657--669.

\bibitem{taylor-law2}
V.~De~La~Pena, P.~Doukhan, Y.~Salhi, {A Dynamic Taylor’s law}, Journal of Applied Probability 59~(2) (2022) 584–607.

\bibitem{ALBRECHER2022}
H.~Albrecher, B.~{García Flores}, Asymptotic analysis of generalized greenwood statistics for very heavy tails, Statistics \& Probability Letters 185 (2022) 109429.

\bibitem{epidemology}
H.~{Scott Hurd}, J.~B. Kaneene, The application of simulation models and systems analysis in epidemiology: a review, Preventive Veterinary Medicine 15~(2) (1993) 81--99.

\bibitem{genetics}
M.~C. Riley, A.~Clare, R.~D. King, Locational distribution of gene functional classes in arabidopsis thaliana, BMC Bioinformatics 8~(1) (2007) 112.

\bibitem{genomics}
B.~D. Peterson-Burch, D.~Nettleton, D.~F. Voytas, Genomic neighborhoods for arabidopsisretrotransposons: a role for targeted integration in the distribution of the metaviridae, Genome Biology 5~(10) (2004) R78.

\bibitem{biology}
S.~L. DeRuiter, I.~L. Boyd, D.~E. Claridge, C.~W. Clark, C.~Gagnon, B.~L. Southall, P.~L. Tyack, Delphinid whistle production and call matching during playback of simulated military sonar, Marine Mammal Science 29~(2) (2013) E46--E59.

\bibitem{economics}
F.~Moscone, E.~Tosetti, Testing for error cross section independence with an application to us health expenditure, Regional Science and Urban Economics 40~(5) (2010) 283--291, advances In Spatial Econometrics.

\bibitem{insurance}
N.~Benlagha, W.~Hemrit, Does investment in insurance stocks reap diversification benefits? static and time varying copula modeling, Communications in Statistics - Simulation and Computation 52~(4) (2023) 1384--1402.

\bibitem{hydrology}
J.~del Castillo, I.~Serra, Likelihood inference for generalized pareto distribution, Computational Statistics \& Data Analysis 83 (2015) 116--128.

\bibitem{optimization}
D.~Vermetten, B.~van Stein, F.~Caraffini, L.~L. Minku, A.~V. Kononova, Bias: A toolbox for benchmarking structural bias in the continuous domain, IEEE Transactions on Evolutionary Computation 26~(6) (2022) 1380--1393.

\bibitem{physics}
P.~Eller, L.~Shtembari, A goodness-of-fit test based on a recursive product of spacings, Journal of Instrumentation 18~(03) (2023) P03048.

\bibitem{materials}
R.~Pakyari, N.~Balakrishnan, Goodness-of-fit tests for progressively type-ii censored data from location–scale distributions, Journal of Statistical Computation and Simulation 83~(1) (2013) 167--178.

\bibitem{anomaly}
O.~Shchur, A.~C. Turkmen, T.~Januschowski, J.~Gasthaus, S.~G\"{u}nnemann, Detecting anomalous event sequences with temporal point processes, in: M.~Ranzato, A.~Beygelzimer, Y.~Dauphin, P.~Liang, J.~W. Vaughan (Eds.), Advances in Neural Information Processing Systems, Vol.~34, Curran Associates, Inc., 2021, pp. 13419--13431.

\bibitem{internet}
C.~Neves, M.~I. Fraga~Alves, Testing extreme value conditions — an overview and recent approaches, REVSTAT-Statistical Journal 6~(1) (2008) 83–100.

\bibitem{athletics}
L.~Henriques-Rodrigues, M.~I. Gomes, D.~Pestana, Statistics of extremes in athletics, REVSTAT-Statistical Journal 9~(2) (2011) 127–153.

\bibitem{emberchts1997}
P.~Emberchts, C.~Kl{\"u}ppelberg, T.~Mikosch, Modelling Extremal Events: for Insurance and Finance, Springer, 1997.

\bibitem{stable1}
J.~P. Nolan, Numerical calculation of stable densities and distribution functions, Communications in Statistics. Stochastic Models 13~(4) (1997) 759--774.

\bibitem{student2}
Student, {The Probable Error of a Mean}, Biometrika 6~(1) (1908) 1--25.

\bibitem{stable3}
G.~Samoradnitsky, M.~S. Taqqu, Stable Non-Gaussian Random Processes: Stochastic Models with Infinite Variance, Routledge, 1994.

\bibitem{shao22}
C.~Nikias, M.~Shao, Signal Processing with Alpha-Stable Distributions and Applications, John Wiley and Sons, New York, 1995.

\bibitem{alek_book}
A.~Janicki, A.~Weron, Simulation and Chaotic Behavior of Alpha-stable Stochastic Processes, Marcel Dekker, Inc., New York, 1994.

\bibitem{non_gauss}
E.~J. Wegman, S.~C. Schwartz, J.~B. Thomas (Eds.), Topics in Non-Gaussian Signal Processing, New York: Springer, 1989.

\bibitem{Nolan2020}
J.~P. Nolan, Univariate Stable Distributions. Models for Heavy Tailed Data, Springer, 2020.

\bibitem{SWtest}
S.~S. Shapiro, M.~B. Wilk, An analysis of variance test for normality (complete samples), Biometrika 52~(3/4) (1965) 591--611.

\bibitem{SWmod1}
D.~A. Pierce, R.~J. Gray, {Testing normality of errors in regression models}, Biometrika 69~(1) (1982) 233--236.

\bibitem{SWmod2}
P.~Royston, Approximating the shapiro-wilk w-test for non-normality, Statistics and Computing 2~(3) (1992) 117--119.

\bibitem{JBtest}
C.~M. Jarque, A.~K. Bera, A test for normality of observations and regression residuals, International Statistical Review / Revue Internationale de Statistique 55~(2) (1987) 163--172.

\bibitem{DPtest}
R.~D'Agostino, E.~S. Pearson, {Tests for Departure from Normality. Empirical Results for the Distributions of $b2$ and $ \sqrt{b1} $}, Biometrika 60~(3) (1973) 613--622.

\bibitem{KStest}
K.~A. Gnedenko~BV, {Limit Distributions of Sums of Independent Random Variables}, Cambridge: Addison-Wesley, 1954.

\bibitem{CMtest}
T.~W. Anderson, {On the Distribution of the Two-Sample Cramer-von Mises Criterion}, The Annals of Mathematical Statistics 33~(3) (1962) 1148 -- 1159.

\bibitem{Kuipertest}
N.~H. Kuiper, Tests concerning random points on a circle, Indagationes Mathematicae (Proceedings) 63 (1960) 38--47.

\bibitem{Watsontest}
G.~S. Watson, Goodness-of-fit tests on a circle, Biometrika 48~(1/2) (1961) 109--114.

\bibitem{ADtest}
T.~W. Anderson, D.~A. Darling, {Asymptotic Theory of Certain "Goodness of Fit" Criteria Based on Stochastic Processes}, The Annals of Mathematical Statistics 23~(2) (1952) 193 -- 212.

\bibitem{LFtest}
H.~W. Lilliefors, On the kolmogorov-smirnov test for normality with mean and variance unknown, Journal of the American Statistical Association 62~(318) (1967) 399--402.

\bibitem{LFmod1}
P.~Sulewski, Modified lilliefors goodness-of-fit test for normality, Communications in Statistics - Simulation and Computation 51~(3) (2022) 1199--1219.

\bibitem{Iskander}
A.~Wyłomańska, D.~R. Iskander, K.~Burnecki, Omnibus test for normality based on the edgeworth expansion, PLOS ONE 15~(6) (2020) 1--36.

\bibitem{NormalityRev2}
K.~R. Das, A.~H. M.~R. Imon, A brief review of tests for normality, American Journal of Theoretical and Applied Statistics 5~(1) (2016) 5--12.

\bibitem{NormalityRev1}
B.~Ebner, N.~Henze, Tests for multivariate normality - a critical review with emphasis on weighted {$L^2$-statistics}, TEST 29~(4) (2020) 845--892.

\bibitem{JelPit2018}
D.~Jelito, M.~Pitera, New fat-tail normality test based on conditional second moments with applications to finance, Statistical Papers 62 (2021) 2083--2108.

\bibitem{infvar1}
L.~Trapani, Testing for (in)finite moments, Journal of Econometrics 191~(1) (2016) 57--68.

\bibitem{infvar2}
I.~Fedotenkov, A bootstrap method to test for the existence of finite moments, Journal of Nonparametric Statistics 25~(2) (2013) 315--322.

\bibitem{ECEM2015}
K.~Burnecki, A.~Wyłomańska, A.~Chechkin, Discriminating between light- and heavy-tailed distributions with limit theorem, PLOS ONE 10~(12) (2015) 1--23.

\bibitem{maraj2023}
K.~Maraj-Zygmat, G.~Sikora, M.~Pitera, A.~Wyłomańska, Goodness-of-fit test for stochastic processes using even empirical moments statistic, Chaos: An Interdisciplinary Journal of Nonlinear Science 33~(1) (2023) 013128.

\bibitem{mssp2023}
K.~Skowronek, T.~Barszcz, J.~Antoni, R.~Zimroz, A.~Wyłomańska, Assessment of background noise properties in time and time–frequency domains in the context of vibration-based local damage detection in real environment, Mechanical Systems and Signal Processing 199 (2023) 110465.

\bibitem{TailidxHill}
B.~M. Hill, A simple general approach to inference about the tail of a distribution, The Annals of Statistics 3~(5) (1975) 1163--1174.

\bibitem{Tailidx2004}
I.~B. Aban, M.~M. Meerschaert, Generalized least-squares estimators for the thickness of heavy tails, Journal of Statistical Planning and Inference 119~(2) (2004) 341--352.

\bibitem{Tailidx2012}
J.~Beran, B.~Das, D.~Schell, On robust tail index estimation for linear long-memory processes, Journal of Time Series Analysis 33~(3) (2012) 406--423.

\bibitem{Tailidx2014}
Y.~M. Tripathi, S.~Kumar, C.~Petropoulos, Improved estimators for parameters of a pareto distribution with a restricted scale, Statistical Methodology 18 (2014) 1--13.

\bibitem{Tailidx2020}
L.~N{\'e}meth, A.~Zempl{\'e}ni, Regression estimator for the tail index, Journal of Statistical Theory and Practice 14~(3) (2020) 48.

\bibitem{Tailidx2023}
J.~Nicolau, P.~M. Rodrigues, M.~Z. Stoykov, Tail index estimation in the presence of covariates: Stock returns’ tail risk dynamics, Journal of Econometrics 235~(2) (2023) 2266--2284.

\bibitem{GenParetoTest1}
J.~Chu, O.~Dickin, S.~Nadarajah, A review of goodness of fit tests for pareto distributions, Journal of Computational and Applied Mathematics 361 (2019) 13--41.

\bibitem{stableTest1}
M.~Pitera, A.~Chechkin, A.~Wy{\l}oma{\'n}ska, Goodness-of-fit test for $\alpha$-stable distribution based on the quantile conditional variance statistics, Statistical Methods \& Applications 31~(2) (2022) 387--424.

\bibitem{stableTest2}
M.-C. Beaulieu, J.-M. Dufour, L.~Khalaf, Exact confidence sets and goodness-of-fit methods for stable distributions, Journal of Econometrics 181~(1) (2014) 3--14, heavy Tails and Stable Paretian Distributions.

\bibitem{csorgo1988}
M.~Cs\"org\H{o}, L.~Horv\'ath, {Asymptotic representations of self-normalized sums}, Probab. Math. Statist. 9~(1) (1988) 15--24.

\bibitem{Hasofer1992}
A.~M. Hasofer, Z.~Wang, A test for extreme value domain of attraction, Journal of the American Statistical Association 87~(417) (1992) 171.

\bibitem{Neves2007}
C.~Neves, M.~I. Fraga~Alves, Semi-parametric approach to the hasofer–wang and greenwood statistics in extremes, TEST 16~(2) (2007) 297–313.

\bibitem{Rivest_1982}
L.~Rivest, Products of random variables and star‐shaped ordering, Canadian Journal of Statistics 10~(3) (1982) 219–223.

\bibitem{ecem2023}
K.~Skowronek, R.~Zimroz, A.~Wyłomańska, Testing for finite variance in time series and its time-frequency representations. applications to vibration signals from rotating machines, Submitted (2023).

\bibitem{student1}
B.~L. Welch, {`Student' and Small Sample Theory}, Journal of the American Statistical Association 53~(284) (1958) 777--788.

\bibitem{Davison_1990}
A.~C. Davison, R.~L. Smith, Models for exceedances over high thresholds, Journal of the Royal Statistical Society: Series B (Methodological) 52~(3) (1990) 393–425.

\end{thebibliography}

\vspace{1cm}

\appendix
\section{Probability distributions}\label{AppA}

In this section we briefly introduce distributions used in testing procedure, namely $\alpha$-stable distribution, Student's t distribution and generalized Pareto distribution.

\subsection{The $\alpha$-stable distribution} \label{App.Stable}

The $\alpha$-stable distribution $\mathcal{S}(\alpha,\beta,\sigma,\mu)$ is defined by its characteristic function (see, e.g., \cite{stable1,stable3})

\begin{equation}
    \Phi(t) = \begin{cases} 
    \text{exp} \left( -\sigma^\alpha |t|^\alpha [1 - i\beta\text{tan} \frac{\pi\alpha}{2} \text{sign}(t) ] + i\mu x\right), & \alpha \neq 1 \\
    \text{exp} \left( -\sigma |t| [1 + i\beta \frac{2}{\pi} \text{sign}(x) \text{log}|t|] + i\mu x \right), & \alpha = 1
    \end{cases}
\end{equation}

\noindent where $\alpha \in (0,2]$ is a stability index responsible for heaviness of the tail of the distribution, $\beta \in [-1,1]$ is a skewness parameter, $\sigma > 0$ is a scale parameter and $\mu \in \mathbb{R}$ is a shift parameter. 

\indent In this work, we analyze the case of symmetric $\alpha$-stable distribution, that is when $\beta = \mu = 0$. In that case the characteristic function reduces to the following formula

\begin{equation}
     \Phi(t) = \text{exp} \left( - \sigma^\alpha |t|^\alpha\right).
\end{equation}
To simplify the notation, we write $\mathcal{S}(\alpha, \sigma)$ instead of $\mathcal{S}(\alpha, 0, \sigma, 0)$.

Note that, for $\alpha \in (0,2)$ tails of the $\alpha$-stable distribution behave as a power-law functions, that is for random variable $X$ with $\mathcal{S}(\alpha, \sigma)$ distribution, we have $\mathbb{P}(X > t) = \mathbb{P}(X < - t) = C t^{-\alpha}(1+ o(1))$, as $t \to \infty$, for some $C > 0$ (see, e.g., Property 1.2.15 in \cite{stable3}). In this case, the variance of the distribution is infinite.
In case of  $\alpha=2$, the $\alpha$-stable distribution reduces to Gaussian distribution $\mathcal{N}(\mu, \sigma^2)$, that is the only case in which the variance exists. \\

\subsection{Student's t distribution} \label{App.Student}

The Student's t distribution $\mathcal{T}(\nu)$ is defined by the following probability density function (see, e.g., \cite{student1,student2})

\begin{equation}
    f(t) = \frac{1}{\sqrt{\nu \pi}} \frac{\Gamma \left( \frac{\nu+1}{2} \right)}{\Gamma \left( \frac{\nu}{2} \right)} {\left(1 + \frac{t^2}{\nu}\right)^{-\frac{\nu+1}{2}}},~~t\in \mathbb{R},
\end{equation}

\noindent where $\nu \in \mathbb{N}_{+}$ represents number of degrees of freedom responsible for the tail behaviour of the distribution. When $\nu \leq 2$, the variance of the distribution does not exist. As $\nu \rightarrow \infty$, the Student's t distribution tends to the Gaussian distribution $\mathcal{N}(0,1)$. In such a case, we use the notation $\mathcal{T}(\infty)$.

\subsection{Generalized Pareto distribution} \label{App.GPD}

Random variable with Generalized Pareto $\mathcal{GP}(\gamma, \delta)$ distribution arise in Peak Over Threshold theory as an approximation for the distribution of a size of exceedances over high threshold. (see, e.g., \cite{Davison_1990} and \cite{emberchts1997}). It is defined through the probability density function

\begin{equation}
    f(t) = \begin{cases}
        \frac{1}{\delta} (1 + \frac{\gamma t}{\delta})^{-1-\frac{1}{\gamma}}, & \gamma \neq 0 \\
    \frac{1}{\delta} \text{exp} (-\frac{t}{\delta}), & \gamma = 0,
    \end{cases}
\end{equation}

\noindent 
where $x \ge 0$ for $\gamma > 0$ and $0 \le x < -\frac{\delta}{\gamma}$ for $\gamma < 0$.
Note, that $\gamma \in \mathbb{R}$ is a parameter responsible for the heaviness of the tail and $\delta > 0$ is a scale parameter. It should be also noted that when $\gamma = 0$, the generalized Pareto distribution reduces to exponential distribution, hence it is a boundary between light tail and heavy tail distributions. The variance of the generalized Pareto  distribution exists if $\gamma < 0.5$.
\section{Quantile tables}
\label{tables}

In this section, we present obtained quantiles of $S_n$ statistic that determine the rejection regions of the proposed tests. The quantiles are obtained using $100000$ Monte Carlo simulations.

\begin{table}[htp!]
\caption{The quantile $Q_{1-c}(n)$ of $S_n$ statistic for Gaussian distribution for different sample sizes and different significance levels $c$.}
\centering
\begin{tabular}{lll}
\hline
n & $c=0.05$  & $c=0.01$ \\ \hline\hline
10 & 0.2125  & 0.2468 \\ \hline
50 & 0.0365  & 0.0387 \\ \hline
100 & 0.0175 & 0.0182 \\ \hline
200 & 0.0085 & 0.0087 \\ \hline
500 & 0.0033 & 0.0033 \\ \hline
1000 & 0.0016 & 0.0016
\end{tabular}
 \label{Tab:gauss}
\end{table}

\begin{table}[htp!]
\caption{The quantile $Q_{c}(n)$ of $S_n$ statistic for generalized Pareto distribution with $\gamma=0.5$ for different sample sizes and different significance levels $c$.}
\centering
\begin{tabular}{lll}
\hline
n & $c=0.05$ & $c=0.01$ \\ \hline \hline
10 & 0.1456 & 0.1300 \\ \hline 
50 & 0.0247 & 0.0220 \\ \hline
100 & 0.0143 & 0.0127 \\ \hline
200 & 0.0068 & 0.0061 \\ \hline
500 & 0.0039 & 0.0035 \\ \hline
1000 & 0.0022 & 0.0021
\end{tabular}
\label{Tab:infvarGP}
\end{table}

\begin{table}[htp!]
\caption{The quantile $Q_{c}(n)$ of $S_n$ statistic for Student's t distribution with $\nu=2$ for different sample sizes and different significance levels $c$.}
\centering
\begin{tabular}{lll}
\hline
n & $c=0.05$ & $c=0.01$ \\ \hline \hline
10 & 0.1312 & 0.1205 \\ \hline
50 & 0.0345 & 0.0316 \\ \hline
100 & 0.0192 & 0.0177 \\ \hline
200 & 0.0106 & 0.0098 \\ \hline
500 & 0.0049 & 0.0045 \\ \hline
1000 & 0.0027 & 0.0025
\end{tabular}
\label{Tab:infvarT}
\end{table}

\section{Proofs}
In this section, we present detailed proofs of Proposition \ref{SnMain}, Proposition \ref{SnOrder}, and Proposition \ref{MGT1.prop}.
\subsection{Proof of Proposition \ref{SnMain}} \label{SnMain.proof}

Note, that for a symmetric random variable $X$ it holds that 
$F_{|X|}(t) = 2F_X(t) - 1$ and hence 
$F^{-1}_{|X|}(u) = F^{-1}_X\left(\frac{1}{2} + \frac{u}{2}\right)$. Therefore,
for any $X^{(\theta_1)} \le_\ast X^{(\theta_2)}$ with distributions $\mathcal{P}_{\theta_1}$ and $\mathcal{P}_{\theta_2}$, respectively, we find that
\[
        \frac{F^{-1}_{|X^{(\theta_1)}|}(u)}{F^{-1}_{|X^{(\theta_2)}|}(u)}
    =
        \frac{F^{-1}_{X^{(\theta_1)}}\left(\frac{1}{2} + \frac{u}{2}\right)}{F^{-1}_{X^{(\theta_2)}}\left(\frac{1}{2} + \frac{u}{2}\right)}
\]
is increasing function with respect to $u$, consequently 
$|X^{(\theta_1)}| \le_\ast |X^{(\theta_2)}|$. 
In order to complete the proof, it suffices to apply Theorem 1 in \cite{arendarczyk2022} 
for random samples $|X_1^{(\theta_1)}|, |X_2^{(\theta_1)}|,\ldots, |X_n^{(\theta_1)}|$
and $|X_1^{(\theta_2)}|, |X_2^{(\theta_2)}|,\ldots, |X_n^{(\theta_2)}|$ with common distributions
$\mathcal{P}_{\theta_1}$ and $\mathcal{P}_{\theta_2}$, respectively.

\subsection{Proof of Proposition \ref{SnOrder}} \label{SnOrder.proof}
For the proof of the case (i) observe that for non-negative random variables $X_1, X_2, \ldots, X_n$, statistics $S_n$ and $T_n$ are equivalent. Thus (i) follows directly from Proposition 1 (i) in
\cite{arendarczyk2022} combined with Theorem 1 in \cite{arendarczyk2022}.

For the case of (ii) of random variables with distribution $S(\alpha, \sigma)$  $\alpha \in (0,2]$, note that it is shown in \cite{Rivest_1982} (see Example 2 in \cite{Rivest_1982}) that for independent random variables $X, Y$, with distributions $S(\alpha_1, \sigma), \alpha_1 \in (0,2)$ and $S(\alpha_2, \sigma), \alpha_2 \in (0,2)$, respectively, inequality $\alpha_1 \le \alpha_2$ implies $Y \le_\ast X$.
For the case $X \sim S(\alpha, \sigma), \alpha \in (0,2)$ and $Y \sim S(2, \sigma)$, that is $Y$ having standard normal distribution, observe that due to proposition 1.3.1 in \cite{stable3} $X$ and $Y$ have the following stochastic representations: $X\stackrel{d}{=} \sqrt{Z}N$ and $Y \stackrel{d}{=} 1\cdot N$, where $Z$ is a positive $\frac{\alpha}{2}$-stable random variable, with $\alpha \in (0,2)$ and $N$ is independent of $Z$, random variable with distribution $\mathcal{N}(0, \sigma^2)$. Moreover, for any random variable $Z$, we have $1 \le_\ast \sqrt{Z}$. Hence, application of Proposition 2 in \cite{Rivest_1982} implies $Y \le_\ast X$, which combined with Theorem \ref{SnMain} completes the proof of the case (iii).

Similar approach can be used  for the case (iii) of Student's t-distribution $\mathcal{T}(\nu)$ with $\nu \in \mathbb{N}$ degrees of freedom. It is shown in \cite{Rivest_1982} (see Example 1 in \cite{Rivest_1982}) that for independent random variables $X, Y$, with distributions $\mathcal{T}({\nu_1})$ and $\mathcal{T}({\nu_2})$, $\nu_1, \nu_2 \in \mathbb{N}$, respectively, inequality $\nu_1 \le \nu_2$ implies $Y \le_\ast X$.
For the case $X \sim \mathcal{T}(\nu), \nu \in \mathbb{N}$ and $Y \sim \mathcal{T}(\infty)$, that is $Y$ having standard normal distribution, note that $X$ and $Y$ can be represented using following stochastic representations: $X\stackrel{d}{=} Z^{-\frac{1}{2}}N$ and $Y \stackrel{d}{=} 1\cdot N$, where $Z$ is a random variable with chi-squared distribution with $\nu$ degrees of freedom and $N$ is independent of $Z$, random variable with standard normal distribution. Moreover, for any random variable $Z$, we have $1 \le_\ast Z^{-\frac{1}{2}}$. Hence, application of Proposition 2 in \cite{Rivest_1982} implies $Y \le_\ast X$. By application of Theorem \ref{SnMain} the proof for the case (ii) is completed.

\subsection{Proof of Proposition \ref{MGT1.prop}} \label{MGT1.prop.proof}
By the definition (\ref{MGT1.rejreg}) of the rejection region for the test $MGT_1$, the power function is of the form
\[
    \mathbb{P}(S_n^{(\alpha)} > \hat{Q}_c(n)),
\]
as a function of $\alpha$, where $S_n^{(\alpha)}$ is given by (\ref{def.sn}) with $X_1, X_2, \ldots, X_n$ heaving a common distribution $\mathcal{S(\alpha, \sigma)}$.
Consequently, decreasing of the power function is a direct consequence of Proposition \ref{SnOrder} (ii). Moreover, by the definition of $\hat{Q}_c(n)$ it is evident that
\begin{eqnarray} \label{size.c}
        \mathbb{P}(S_n^{(2)} > \hat{Q}_c(n)) = c.
\end{eqnarray}
Hence, the test $MGT_1$ has size c. 
In order to prove that the test $MGT_1$ is unbiased it suffices to note that stochastic decreasing of the power function in conjuction with (\ref{size.c}) implies that 
\[
        \mathbb{P}(S_n^{(\alpha)} > \hat{Q}_c(n))
    \ge 
        c
\] 
for all $\alpha \in (0,2)$. This completes the proof.

\end{document}